\documentclass{amsart}
\usepackage{amssymb,amsmath,amscd,xy,graphicx,textcomp}
\newtheorem{theorem}{Theorem}[section]
\newtheorem{lemma}[theorem]{Lemma}

\newtheorem{definition}[theorem]{Definition}

\newtheorem{proposition}[theorem]{Proposition}

\xyoption{arrow}

\xyoption{matrix}

\def\C{\mathbb{C}}

\def\Z{\mathbb{Z}}

\title{Bounds on generators and relations for the algebra of $SL_2(\C)$ conformal blocks}
\author{Christopher Manon}
\thanks{}

\begin{document}

\begin{abstract}
We show that the Cox ring of the moduli of $SL_2(\C)$ quasi-parabolic principal bundles on a marked curve 
is generated by conformal blocks of level $1$ and $2.$  We show that the ideal which vanishes on these
generators is generated by forms of degrees $2, 3, 4.$ 
\end{abstract}

\maketitle

\tableofcontents

\smallskip

\section{Introduction}

In this note we give a combinatorial argument which bounds the degree of generators and relations needed to present the Cox ring
of the moduli of $SL_2(\C)$ quasi-parabolic principal bundles on a smooth, projective curve with marked points.  It should be read 
as a continuation of the author's work in \cite{M4}.

Let $(C, \vec{p}) \in \mathcal{M}_{g, n}$ be a smooth, projective curve with $n$ chosen marked points $\vec{p} \subset C.$
We study the structure of the Cox ring (or total coordinate ring) of the moduli stack $\mathcal{M}_{C, \vec{p}}(SL_2(\C))$ of quasi-parabolic $SL_2(\C)$
principal bundles on $(C, \vec{p}).$  The Cox ring of a scheme, algebraic space, or stack $M$ can be roughly described as the direct sum of the spaces of global sections of all line bundles, with the mutliplication operation given by composition of global sections. 

\begin{equation}
Cox(M) = \bigoplus_{L \in Pic(M)} H^0(M, L)\\
\end{equation} 

 Lazlo, Sorger  have determined that the structure of the Picard group of the stack $\mathcal{M}_{C, \vec{p}}(SL_2(\C))$ is isomorphic
to $\Z^n \times \Z$, \cite{LS}. For a line bundle $\mathcal{L}(\vec{r}, L) \in Pic(\mathcal{M}_{C, \vec{p}}(SL_2(\C))$  to be effective, it is
necessary for $r_i \leq L$. We determine a finite generating set for $Cox(\mathcal{M}_{C, \vec{p}}(SL_2(\C))$, and a finite set of
relations which hold among these generators, when the curve $(C, \vec{p})$ is generic. 

\begin{theorem}\label{main0}
The Cox ring $Cox(\mathcal{M}_{C, \vec{p}}(SL_2(\C)))$ is generated by the sections of line bundles with $L \leq 2$,
and the relations which hold among these generators are generated by those with total $L$ value $2, 3, 4.$
\end{theorem}

Finite generation of the Cox ring is a critical property in the birational theory of a space, in this case
the coarse moduli associated to $\mathcal{M}_{C, \vec{p}}(SL_2(\C))$.  It implies that the space
in question is a Mori Dream space.  This implies that the Mori minimal model program can be carried
out for the space, see \cite{HK}.

\begin{figure}[htbp]
\centering
\includegraphics[scale = 0.65]{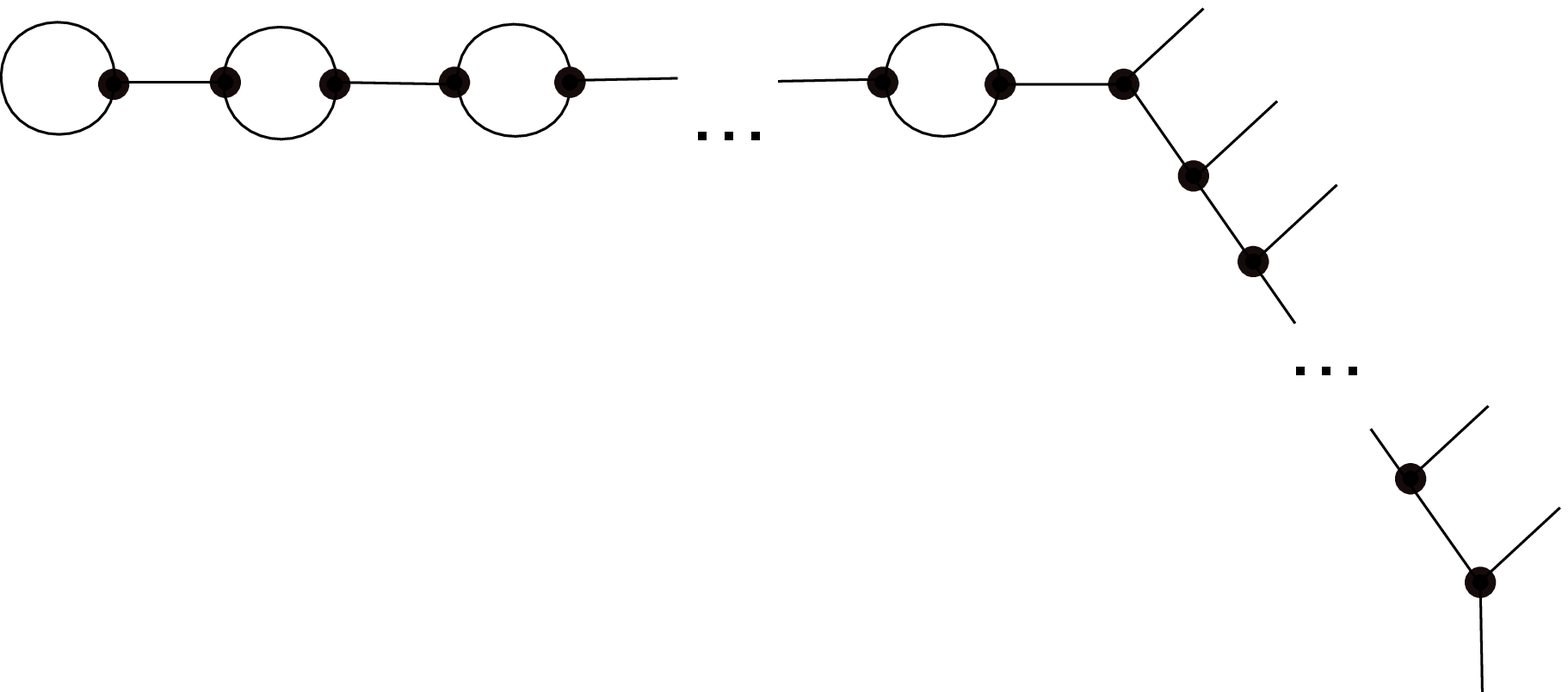}
\caption{}
\label{virus}
\end{figure}

The global sections of the line bundle $\mathcal{L}(\vec{r}, L)$, $(\vec{r}, L) \in \Z^n \times \Z$ can 
be identified with the level $L$ $sl_2(\C)$ conformal blocks with weight data $\vec{r}$ on the curve $(C, \vec{p})$.
Conformal blocks are the correlation functions from the Wess-Zumino-Witten model of conformal field theory. Spaces
of conformal blocks constitute a Topological Quantum Field Theory, which forces them to satisfy certain
combinatorial identities called the factorization rules, see \cite{B}.  In \cite{M4} we used  these rules
to find one flat toric degeneration of $Cox(\mathcal{M}_{C, \vec{p}}(SL_2(\C))$ for every trivalent graph $\Gamma$
with $n$ leaves and first Betti number equal to the genus of $C.$  In order to prove Theorem \ref{main0} we study one particular degeneration.
Let $\Gamma(g, n)$ be the graph depicted in Figure \ref{virus} with $n$ leaves and $g$ loops. 

We consider the graded semigroup defined by the following polytopes.

\begin{definition}
Let $P_{g, n}(L)$ be the polytope defined by non-negative real number weightings of the edges $e \in E(\Gamma(g, n))$
such that the three weights $w_1(v), w_2(v), w_3(v)$ trivalent vertex $v \in V(\Gamma(g, n))$ satisfy the following inequalities.  

\begin{enumerate}
\item $w_1(v) + w_2(v) + w_3(v) \leq 2L$\\ 
\item $|w_1(v) - w_2(v)| \leq w_3(v) \leq w_1(e) + w_2(e)$\\
\end{enumerate}

\end{definition} 

Although it may not be immediately obvious, the reader can verify that these inequalities are symmetric in the entries. We consider
these polytopes with respect to the lattice defined by the conditions that the weights on all edges be integers, and $w_1(v) + w_2(v) + w_3(v) \in 2\Z.$
From now on we abuse notation and refer to the set of lattice points in the polytope by $P_{g, n}(L)$. There is a natural multiplication operation $\pi: P_{g, n}(L) \times P_{g, n}(K) \to P_{g, n}(L+K)$ defined by adding weightings edge-wise. 

\begin{figure}[htbp]
\centering
\includegraphics[scale = 0.65]{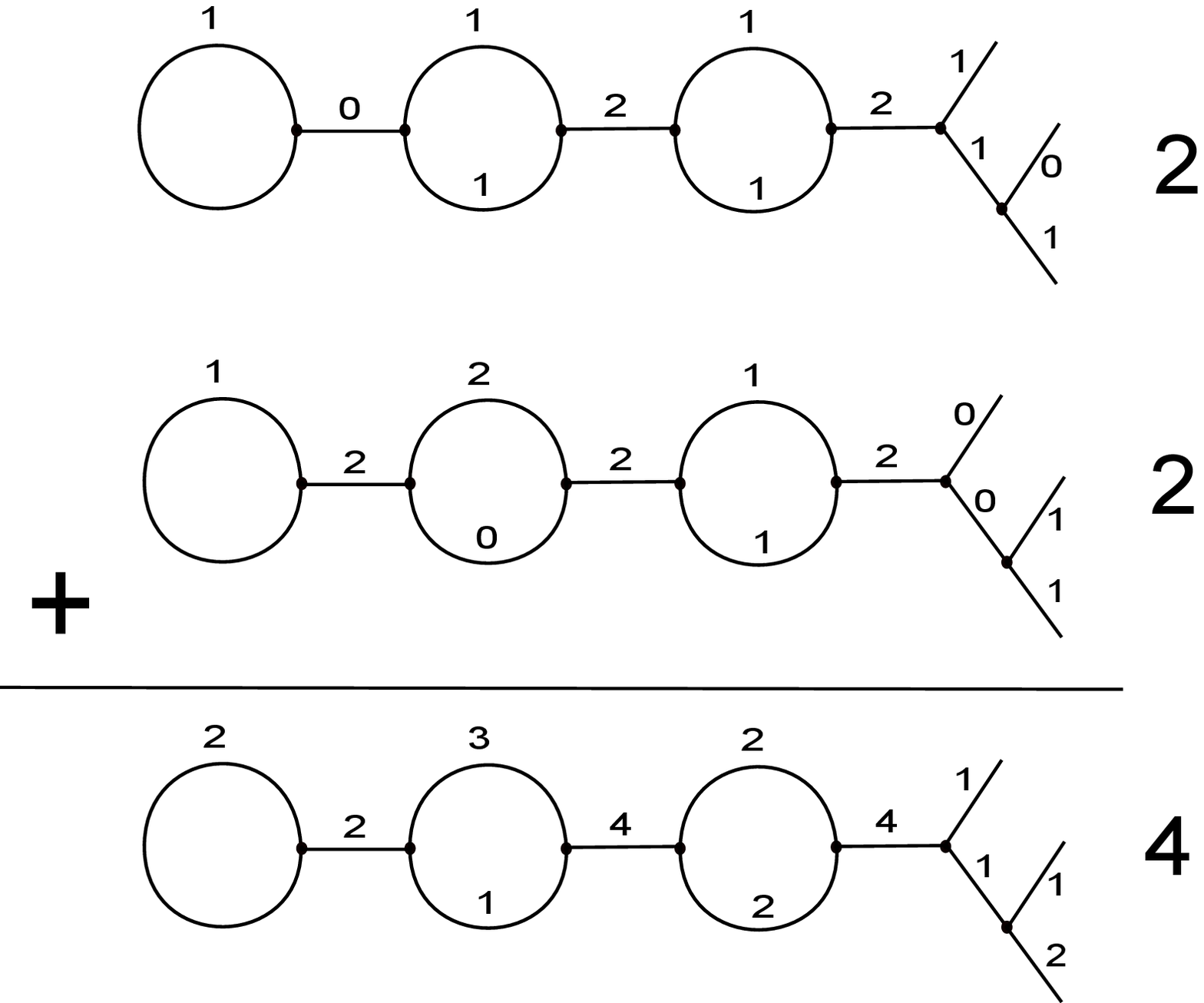}
\caption{}
\label{Fig1}
\end{figure}

We let $P_{g, n}$ be the graded semigroup defined by this multiplication operation, with $\C[P_{g, n}]$ the associated graded semigroup algebra.   
The following is a consequence of the work in \cite{M4}.

\begin{theorem}
Let $(C, \vec{p})$ be a smooth, projective curve of genus $g$ with $n$ marked points, then there is a flat degeneration of algebras,

\begin{equation}
Cox(\mathcal{M}_{C, \vec{p}}(SL_2(\C))) \Rightarrow \C[P_{g, n}]\\
\end{equation}

\end{theorem}

The purpose of this note is to prove the following theorem on the semigroup algebra $\C[P_{g, n}]$. Theorem
\ref{main0} then follows by standard flat degeneration arguments. 

\begin{theorem}\label{main}
The semigroup algebra $\C[P_{g, n}]$ is generated by lattice points in $P_{g, n}(1)$ and $P_{g, n}(2)$, and the ideal
of relations on these generators is generated in degrees $\leq 4.$
\end{theorem}
 
Theorem \ref{main} should be compared to recent results of Buczynska, Buczynski, Kubjas, and Michalek \cite{BBKM}, who 
prove that the relevant semigroup is generated in degree $\leq g+1$ for $any$ graph with $n$ leaves and first Betti number equal to $g.$
This result also extends a theorem of Abe \cite{A}, which establishes the generation statement when $n =0$.

\section{Generators}

We first show that elements from $P_{g, n}(1)$ and $P_{g, n}(2)$ suffice to generate the semigroup $P_{g, n}.$  This
part of Theorem \ref{main} is broken into three lemmas. 

\begin{lemma}\label{g1}
Theorem \ref{main} for case $(g, n)$ follows from cases $(g, 1)$ and $(0, n + 1)$
\end{lemma}

\begin{lemma}\label{2L+1}
The multiplication map $\pi: P_{g, 1}(1)\times P_{g, 1}(2L) \to P_{g, 1}(2L+ 1)$ is surjective.
\end{lemma}

\begin{lemma}\label{2L}
The multiplication map $\pi: P_{g, 1}(2)^{\times L} \to P_{g, 1}(2L)$ is surjective.
\end{lemma}


\subsection{Proof of Lemma \ref{g1}}

We prove Lemma \ref{g1} with a fiber-product argument. Case $(0, n)$ of this theorem is a byproduct of the work in \cite{BW}, where Buczynksa and Wiesniewski show that $P_{0, n}$ is generated in degree $1$ with quadratic relations.  We state the relevant aspects of their theorem below. 

\begin{proposition}
$P_{0, n}$ is generated by weightings from $P_{0, n}(1)$, all of which have a $1$ or a $0$ on a given edge of $\Gamma(0, n)$. 
\end{proposition}  

\begin{figure}[htbp]
\centering
\includegraphics[scale = 0.6]{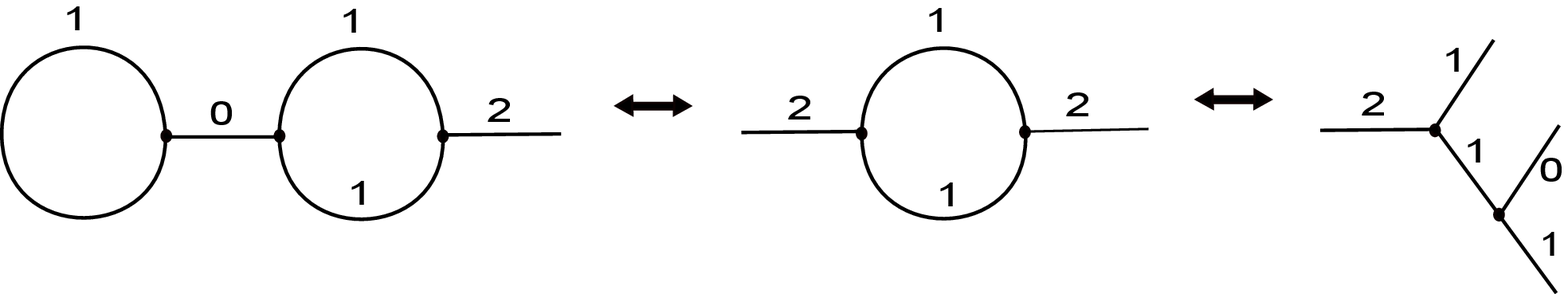}
\caption{}
\label{Fig2}
\end{figure}

We also make use of the fact that the parity of the integers which weight the edges of $\Gamma(g, 1)$ in an element of $P_{g, 1}$ 
are combinatorially restricted.  First, the weight on any edge of $\Gamma(g, n)$ coming from an element of $P_{g, n}(L)$ must be less
than or equal to $L$, this follows directly from the defining inequalities. Second, note that the even-sum condition which defines the lattice
under consideration forces each horizontal edge in $\Gamma(g, 1)$ to be weighted with an even number, including the edge on the far right. 
With these facts in mind, we consider an element of $P_{g, n}$ ($g > 0$) as an element of $P_{g, 1}$ glued to an element of $P_{0, n+1},$
see Figure \ref{Fig2}. 

Let $\omega \in P_{g, n}(L),$  and let $\omega^g \in P_{g, 1}$ and $\omega^0 \in P_{0, n+1}$ be the restrictions
of $\omega$ to the copies of $\Gamma(g, 1)$ and $\Gamma(0, n+1)$ in $\Gamma(g, n)$ respectively.  By the theorem
of Buczynska and Wiesniewski, $\omega^0$ decomposes into $L$ weightings of $\Gamma(0, n+1)$.

\begin{equation}
\omega^0 = \eta_1 + \ldots + \eta_L\\
\end{equation}

Now suppose that the generation statement of Theorem \ref{main} holds for $P_{g, 1},$ then 
$\omega^g$ likewise decomposes into elements of degree $1$ and $2$.

\begin{equation}
\omega^g = [\alpha_1 + \ldots + \alpha_k] + [\beta_1 + \ldots + \beta_m]\\
\end{equation}

\noindent
Here $2m + k = L$.  Note that none of the $\alpha_i$ can have a non-zero weight on the edge which is shared
by $\Gamma(g, 1)$ and $\Gamma(0, n+1)$ in $\Gamma(g, n).$  We let $\eta_1, \ldots, \eta_{L'}$ and $\beta_1, \ldots, \beta_m'$ be the elements
with a non-zero weight on this edge. Then we must have 
$L' = 2m'.$   We simply pair up the $m'$ pairs of elements $\eta_1 + \eta_2, \ldots$ to make $m'$ elements in $P_{0, n+1}(2),$ 
and glue these to the $\beta_i$ along the shared edge.  What remains must be exactly $k'$ $\eta_i'$s and $k'$ $\alpha_i'$s, both with
$0$ along the shared edge, which can then likewise be glued.  This proves Lemma \ref{g1} in the case $n > 0.$  For the case $n= 0$, we simply note that $P_{g, 0} \subset P_{g, 1}$ is the set of weightings with a $0$ on the unique leaf edge, and a factorization of such an element
must result in elements which likewise weight the leaf edge with a $0.$

\subsection{Proof of Lemma \ref{2L+1}}

\begin{figure}[htbp]
\centering
\includegraphics[scale = 0.85]{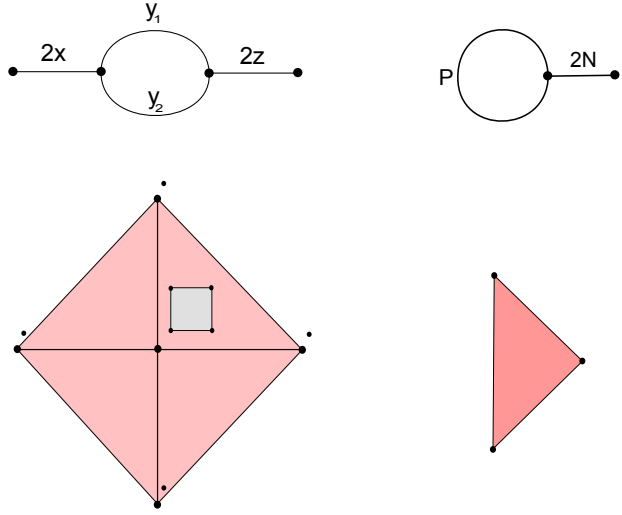}
\caption{$B_1(L), B_2(L)$}
\label{Fig3}
\end{figure}

We prove Lemma \ref{2L+1} with an analysis of weightings on the component graphs in Figure \ref{Fig3}.

\begin{figure}[htbp]
\centering
\includegraphics[scale = 0.65]{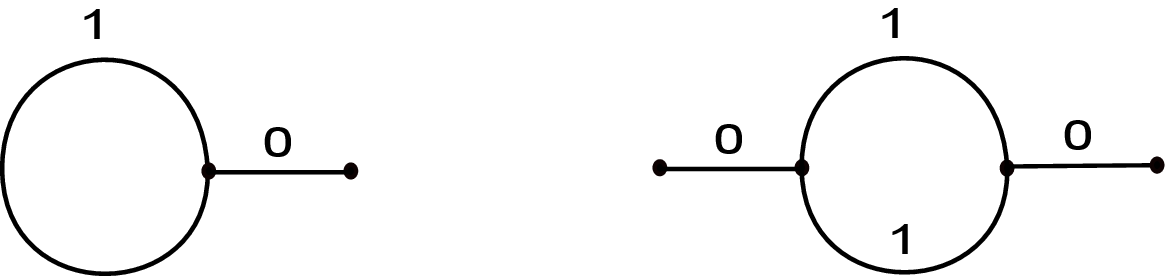}
\caption{}
\label{Fig4}
\end{figure}

We define the graded semigroups $B_1$ and $B_2$ accordingly.  
We consider the weightings of these graphs by single loops $o_1$ and $o_2$, depicted in Figure \ref{Fig4}.
These are the only non-trivial elements in $B_1(1)$ and $B_2(1),$ respectively.  The following proposition
is the crux of both Lemma \ref{2L+1} and \ref{2L}. 

\begin{proposition}\label{bb}
The semigroups $B_1$ and $B_2$ are generated by elements of degrees $1,$ and $2.$
\end{proposition}

\noindent
Note that it follows that all generators weight the leaf edges of either graph with either $0$ or $2.$
We refer to Figure \ref{Fig3} in what follows. 

\begin{figure}[htbp]
\centering
\includegraphics[scale = 0.65]{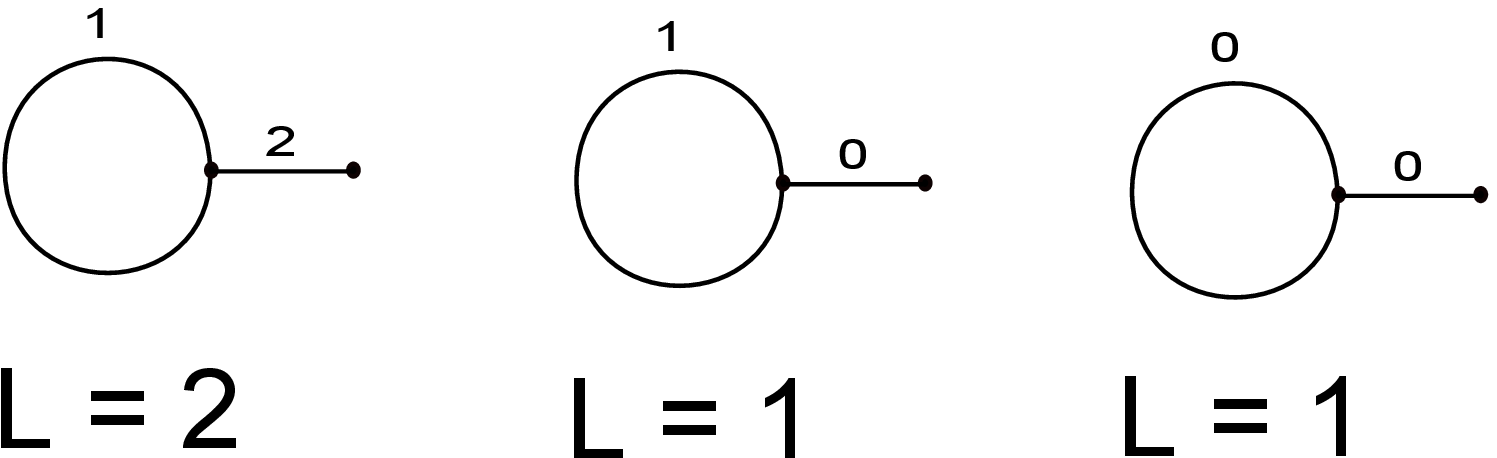}
\caption{}
\label{Fig5}
\end{figure}

\begin{proof}[Proof for $B_1$]
Let $\omega$ be an element of $B_1.$  If $2P > 2N$ then the element $o_1$ can be factored off
of $\omega$.  If $2P = 2N$ then $\omega$ is a power of the element with $P = 1, N = 2,$ with a product
of the element with $P = N = 0,$ depicted
in Figure \ref{Fig5}.
\end{proof}

\begin{figure}[htbp]
\centering
\includegraphics[scale = 0.40]{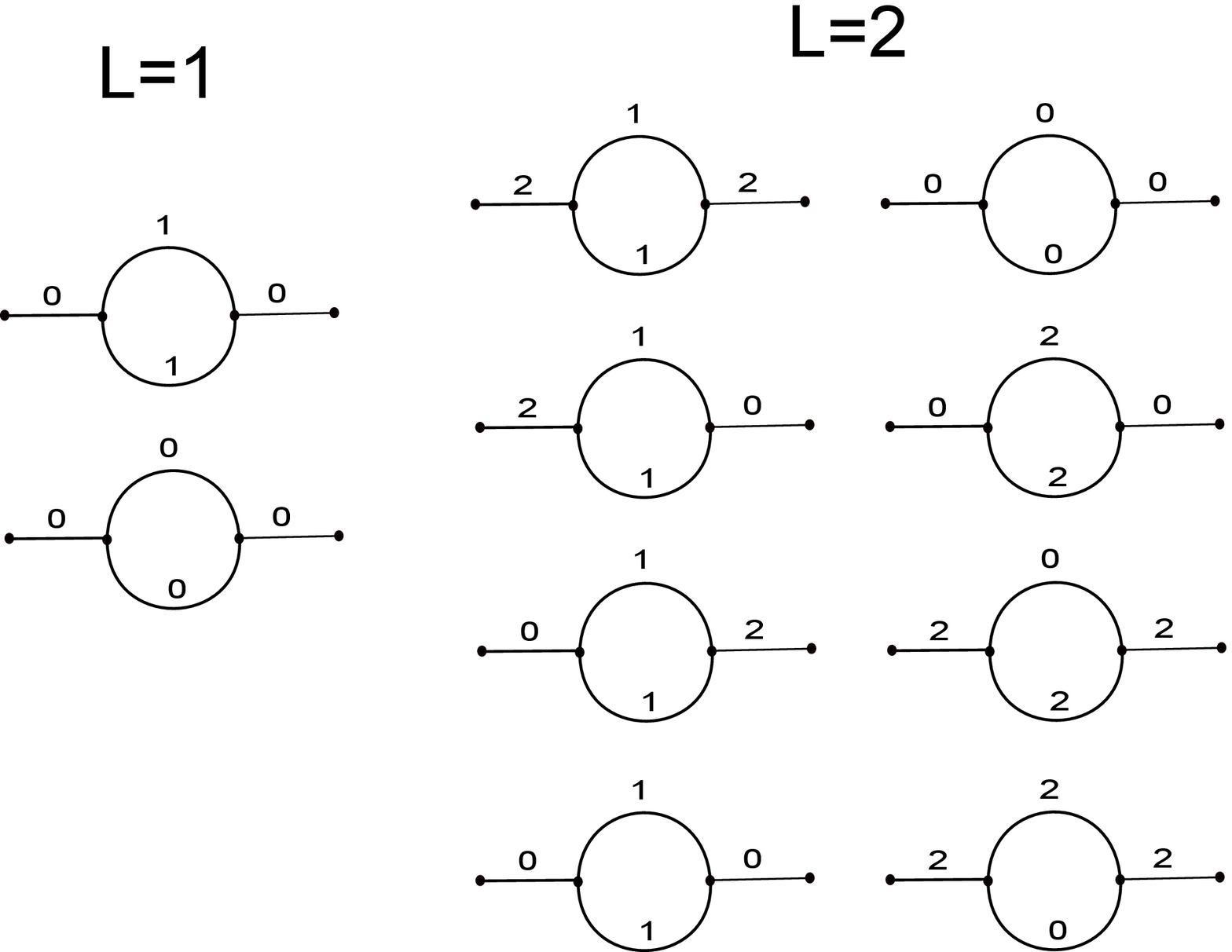}
\caption{}
\label{Fig6}
\end{figure}

\begin{proof}[Proof for $B_2$]
We assume without loss of generality that $y_1 \geq y_2$.  Note that if $y_1 > y_2$
then $y_1 -  y_2 \in 2\Z,$  and by the triangle inequalities, $2x, 2z, \geq 2.$
We can correspondingly factor off an element with $2x = 2z = y_1 = 2.$
This can be done until $y_1 = y_2.$ In this case, if $x > z,$ we can factor
off an element with $ 2x = 2, z = 0, y_1 = y_2 = 1.$ The remainder is a product
of elements with $y_1 = y_2, x = z$. These can always be factored into 
generators with these properties, depicted in Figure \ref{Fig6}. 
\end{proof}

\begin{figure}[htbp]
\centering
\includegraphics[scale = 0.65]{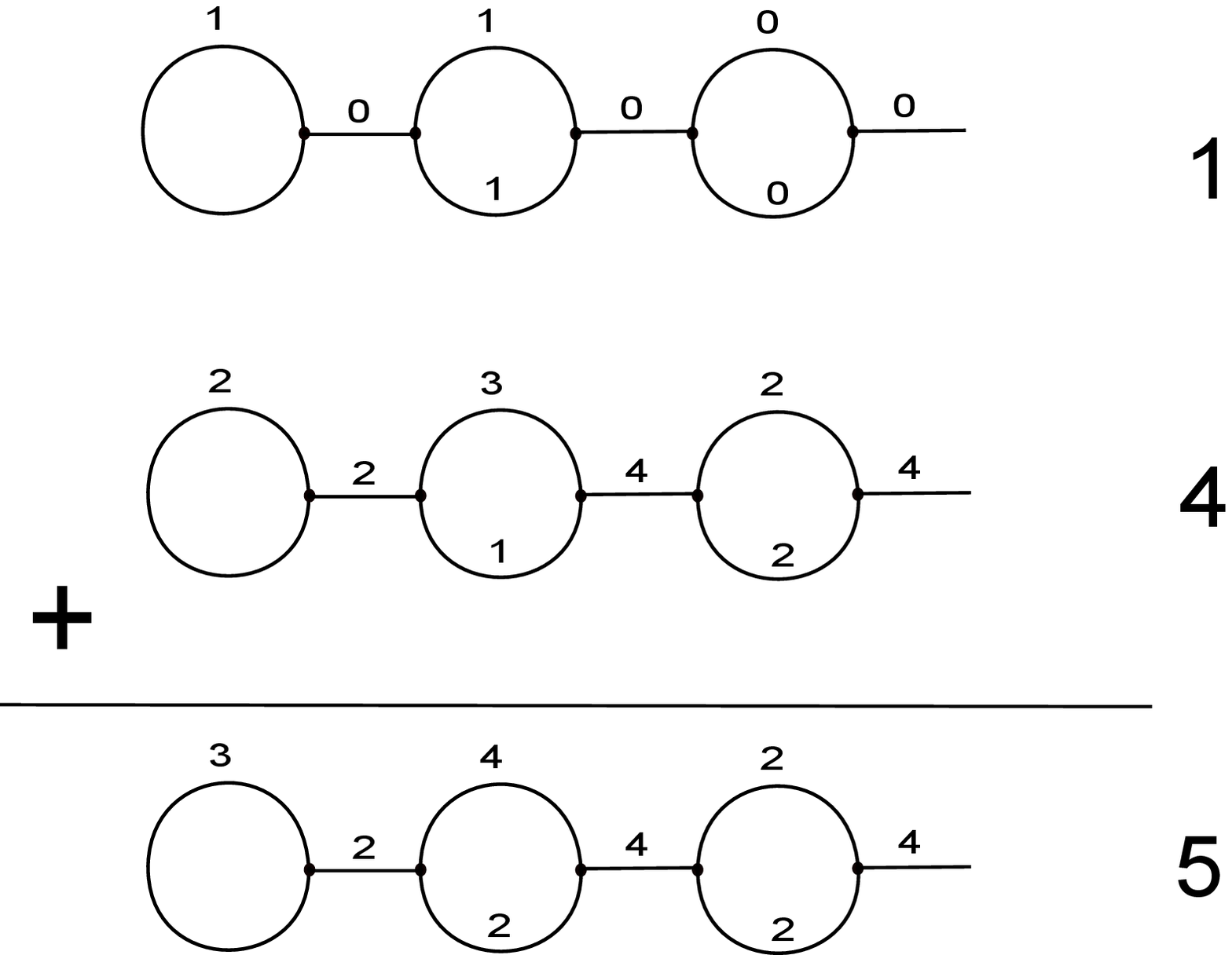}
\caption{}
\label{level1}
\end{figure}

Both of these arguments are captured by results in \cite{M7}. Now if, $\omega \in P_{g, 1}(2L+1),$
then we can check each loop in $\Gamma(g, n)$.  If the loop has a vertex with entries adding to $2L+1$, 
we can extract an $o_1$ or $o_2$ from that loop. If not, we leave it alone.  The resulting element is 
a member of $P_{g, n}(1),$ and the remainder is an element of $P_{g, n}(2L),$ see Figure \ref{level1}.

\subsection{Proof of Lemma \ref{2L}}

To prove Lemma \ref{2L} we note that a weighting of $\Gamma(g, 1)$ automatically
induces a weighting of $\Gamma(g-1, 1)$ and $B_2$.  Now we are in a position to use an inductive argument similar to the proof of Lemma \ref{g1}. 
The theorem is true for $\Gamma(1, 1) = B_1$ by the proof of Proposition \ref{bb}.  Suppose
it holds for $\Gamma(g -1, 1).$  For an element $\omega \in P_{g, 1}(2L)$ we consider
the restrictions $\omega_1, \omega_2$ to $\Gamma(g-1, 1)$ and $B_2$ respectively. 
By induction, $\omega_1$ can be factored into $L$ elements of $P_{g-1, 1}(2).$

\begin{equation}
\omega_1 = \alpha_1 + \ldots + \alpha_L\\
\end{equation}

Similarly, the same thing can be carried out for $\omega_2$. 

\begin{equation}
\omega_2 = \beta_1 + \ldots + \beta_L\\
\end{equation}

The number of $\alpha_i$ which weight the edge shared by $\Gamma(g-1, 1)$ and $B_2$ with a $2$
must then equal the number of $\beta_i$ which also weight this edge $2.$ These elements can be matched, 
leaving only elements which weight the shared edge $0$.  This proves Lemma \ref{2L}.

\section{Relations}

We analyze relations in $P_{g, n}$ in a similar manner
to generators.  First, we will describe relations for the building
block polytopes $B_1, B_2$, then we show that these relation results
are stable as these building blocks are glued together.  This will give
the result for $P_{g, 1}$.  Then we will show that this extends to $P_{g, n}.$

\subsection{The polytopes $B_1$ and $B_2$}

The semigroup algebra $\C[B_1]$ is generated by two
elements of level $1$ and an element of degree $2,$
and no relations hold between these elements.

The semigroup algebra $\C[B_2(2)]$ of even level
elements in $B_2$ was studied in \cite{M7}. It was
shown to have relations generated in degree $2,$ which
translates to level $4.$ To see why this is the case, we can
describe a unimodular triangulation of this polytope directly. 

Following \cite{M7}, we rewrite the elements of $B_2(2),$ taking a weighting ($2a, x, y, 2b)$
to $(a, \frac{1}{2}(x-y), \frac{1}{2}(x+y) - 1, b).$  This results in a polytope $B_2'(2)$ with the following
$8$ lattice points. 

$$[0000] \ \ [ 0001] \ \ [1000] \ \ [1001]$$

$$[0010] \ \ [1101] \ \ [0010] \ \ [1101]$$

\noindent
This can be split into $4$ quadrants, given as the convex hulls of the following lists of points.

$$[0000] \ \ [ 0001] \ \ [1000] \ \ [1001]  \ \ [0010] \ \ [1101]$$

$$[0000] \ \ [ 0001] \ \ [1000] \ \ [1001] \ \ [1101] \ \ [0010]$$

$$[0000] \ \ [ 0001] \ \ [1000] \ \ [1001]  \ \ [1101] \ \ [0010]$$

$$[0000] \ \ [ 0001] \ \ [1000] \ \ [1001]  \ \ [0010] \ \ [1101]$$

These quadrants can be recovered as the regular regions
of the term order given by weighting the lattice point $[A B C D ]$ with $A^2 + B^2 + C^2 + D^2.$
The fact that each quadrant is a normal polytope implies that $B_2(2)$ is "balanced", see \cite{M1}, \cite{M7}.
This means that the degree of necessary relations for $B_2(2)$ is bounded by those of each quadrant.
Each quadrant's relations are generated by a single binomial, $[0001][1000] = [0000][1001]$. 
Now let $\omega \in B_2$ be an element of even level, and let
this element be factored in two distinct ways. 

\begin{equation}
\omega = \sum \beta_i = \sum \alpha_j\\
\end{equation}

Since $B_2$ is generated by elements of degrees $1$ and $2,$
we can collect both sides of this equation into elements
of degree $2$, 

\begin{equation}
\sum \beta_i^* = \sum \alpha_j^*.\\
\end{equation}

\noindent
The expression $\sum \beta_i^*$ 
can be transformed to $\sum \alpha_j^*$ using level $4$ binomials, as above.  
This proves that level $4$ relations suffice in the case of an even
level element. Now suppose $\omega$ has an odd level.  We consider two factorizations
into level $1$ and $2$ elements. 

\begin{equation}
\omega = \sum \beta_i = \sum \alpha_j\\
\end{equation}

\noindent
Note that both factorizations must include at least one level $1$ element. 
If the same level $1$ element appears in both factorizations, it can be pulled
off, and we can apply the argument of the even case, above. Otherwise, we have
two cases.  

We consider the sums of weights around the two internal vertices of 
the graph.  If one of these sums is $2(2L+1),$ each degree $2$ element must 
have this sum equal to $4$, and each degree $1$ element must have this sum equal to $2$, 
in both factorizations.  This implies both factorizations contain the loop
element of level $1$ from Figure \ref{Fig4}, as this is the only element of level $1$ with this property. 

If both sums are strictly less than $2(2L+1)$, then some element on both sides
has these sums less than their maximum possible value ($2$ for degree $1$ elements, $4$ for degree $2$). 
If this element is level $1$, then it must be the element which weights all edges $0$.
If it is level $2$, it can be factored into the element which weights all edges $0$ and the loop element. 
Either way, the empty element can be factored off. This proves the relations component of Theorem  \ref{main} for $B_1$ and $B_2.$

\subsection{Relations for $P_{g, 1}$}

We describe relations for $P_{g, 1}$ by building this semigroup out of copies of $B_1$ and $B_2$
with toric fiber products.   The theorem holds for $B_1$ because of the unique factorization property,
this is the base case of our induction.  To treat $P_{g, 1}$, we consider it as a toric fiber product
of $P_{g-1, 1}$ with $B_2$ over the semigroup of non-negative integer points $(n, L)$ with $2n \leq L$. This semigroup
is generated by $(0, 1)$ and $(1, 2)$, and also has unique factorization.  The maps from $P_{g-1, 1}$ and $B_2$ to this semigroup are computed by returning $\frac{1}{2}$ the weight on the edge shared by the supporting graphs of these semigroups, and the level of the weightings. 
Given an element $\omega \in P_{g, 1}$, and two factorizations

\begin{equation}
\omega = \sum \beta_i = \sum \alpha_j\\
\end{equation}

\noindent
we first collect level $1$ terms into level $2$ terms, to obtain factorizations
into all level $2$ elements in the even level case, and all level $2$ elements
except one level $1$ element in the odd case.  Note that in doing so, we have 
used at most level $2$ relations.

\begin{equation}
\omega = \sum \beta_i^* = \sum \alpha_j^*\\
\end{equation}

Now we consider the restriction of this relation to $P_{g-1, 1}$.  By induction, 
there is a way to transform $\sum \beta_i^*|_{g-1, 1}$ into $ \sum \alpha_j^*|_{g-1, 1}$
using degree $2, 3, 4$ relations.  We claim that each of these relations can be lifted
to a relation in $P_{g, 1}.$  This holds because any such relation must preserve the list of values
along the shared edge.  This means that we can "unglue" the $B_2$ side of the weightings involved
in the relation, perform the relation on the $(g-1, 1)$ side, and glue the $B_2$ weightings back on in 
any way consistent with their values along the shared edge to obtain a relation on $P_{g,1}.$

\begin{figure}[htbp]
\centering
\includegraphics[scale = 0.55]{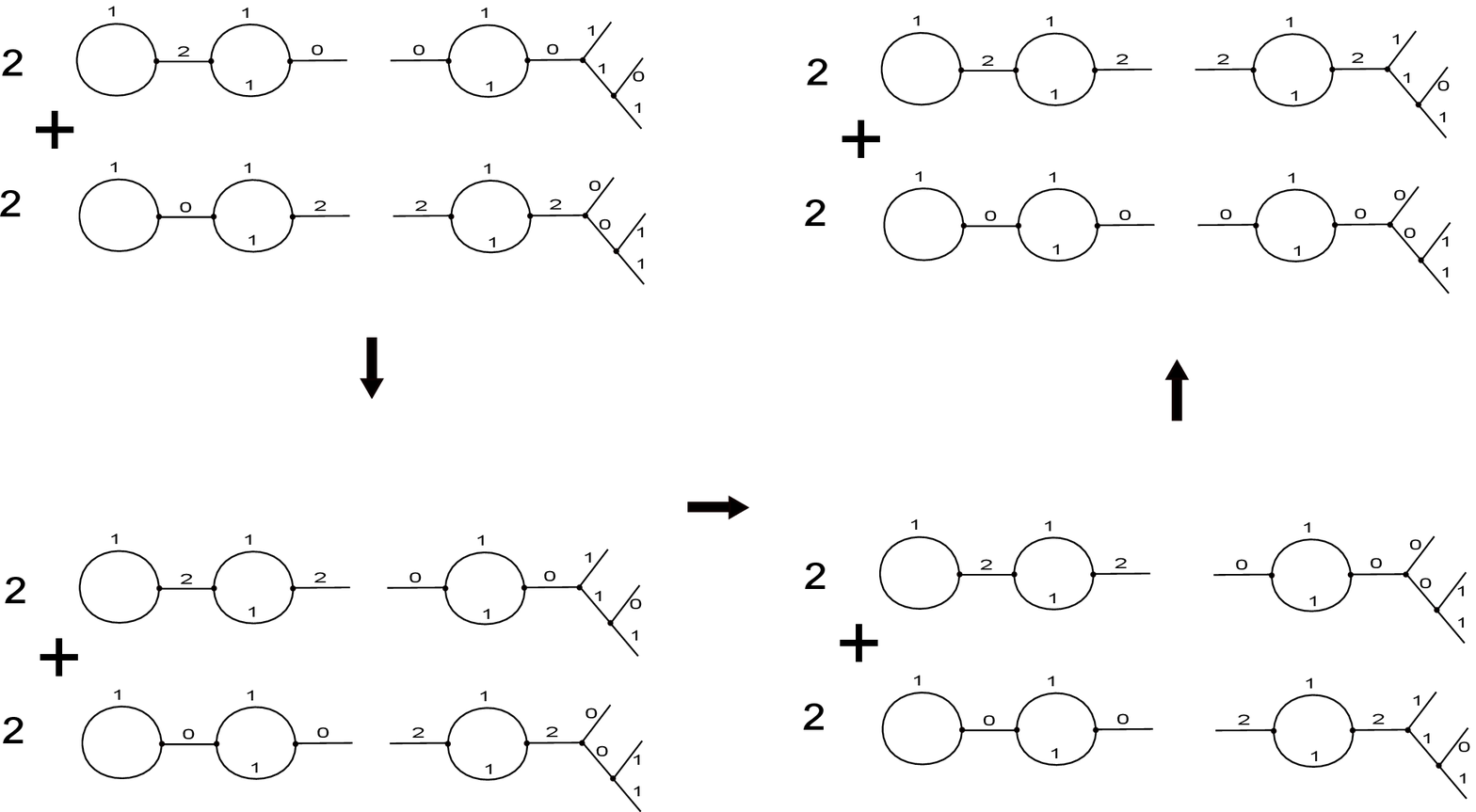}
\caption{}
\label{Fig8}
\end{figure}

This proves that we can transform $\sum \beta_i^*$ to  $\sum \alpha_j^*$ "along $P_{g-1, 1}$."
But we can now apply the same argument to perform this transformation over $B_2$, using the
relations above.  This completes the proof of the relation statement in Theorem \ref{main} 
for $(g, 1),$ and an identical argument applies to $(g, 0).$

\subsection{Relations for $P_{g, n}$}

The same argument used in the case $P_{g, 1}$ to extend relations from $P_{g-1, 1}$ to $P_{g, 1}$ can 
be used to show that any relation on $P_{g, 1}$ can be extended to a relation on $P_{g, n}$.  From this
it follows that given two normalized ways to represent an element $\omega = \sum \beta_i^* = \sum \alpha_j^*,$
one can be transformed to the other "over $P_{g, 1}$" with relations of degree $2, 3, 4.$  Now the same argument
can be applied over the edge which separates $P_{0, n+1}$ from $P_{g, 1}$ using the degree $2$  $P_{0, n+1}$ relations discovered in \cite{BW}.  This completes Theorem \ref{main}.

\bigskip
\noindent
Christopher Manon:\\
cmanon@gmu.edu\\
Department of Mathematics,\\ 
George Mason University,\\ 
Fairfax, VA  22030 USA

\date{\today}

\end{document}